\documentclass[11pt,twoside,a4paper]{article}
\usepackage{graphicx}
\usepackage{amsfonts}
\usepackage{bbm}
\usepackage{mathrsfs}
\usepackage{fancyhdr}
 \usepackage{mathrsfs,amsmath,amssymb,amsthm}
\usepackage{color}
 \usepackage{amssymb}
 \usepackage{amsbsy}

\allowdisplaybreaks

 \theoremstyle{plain}
\theoremstyle{remark}  \newtheorem{remark}{\noindent\mbox{Remark}}
 \theoremstyle{plain}
 \theoremstyle{plain}\newtheorem{lemma}{\noindent\mbox{Lemma}}
\theoremstyle{plain} \newtheorem{theorem}{\noindent\mbox{Theorem}}
 \theoremstyle{plain}\newtheorem{proposition}{\noindent\mbox{Proposition}}
 \theoremstyle{plain}\newtheorem{corollary}{\noindent\mbox{Corollary}}
\theoremstyle{definition} 

 \def\proof{\noindent{\it Proof.~~}}
 \def\ve{\varepsilon}
 \def\qed{\hfill$\Box$\medskip}
 \def\rto{\rightarrow\infty}
\def\z{\left}
\def\y{\right}
 \def\no{\nonumber}

\begin{document}
 \title{\textbf{Times of a branching process with immigration in varying environment attaining a fixed level}\thanks{Supported by National
Natural Science Foundation of China (Grant No. 11501008)}}                  

\author{\uppercase{Hua-Ming Wang}$^{\dag}$
}
\maketitle%
 \footnotetext[2]{Email: hmking@ahnu.edu.cn}
\vspace{-.5cm}

\begin{center}
\begin{minipage}[c]{12cm}
\begin{center}\textbf{Abstract}\quad \end{center}
Consider a branching process $\{Z_n\}_{n\ge 0}$ with immigration in varying environment. For $a\in\{0,1,2,...\},$ let $C=\{n\ge0:Z_n=a\}$ be the collection of times
    at which the population size of the process attains level $a.$ We give a criterion to determine whether  the set $C$ is finite or not.  For critical Galton-Watson process, we show that $|C\cap [1,n]|/\log n\rightarrow S$ in distribution, where $S$ is an exponentially distributed random variable with $P(S>t)=e^{-t},\ t>0.$

\vspace{0.2cm}

\vspace{0.2cm}

\end{minipage}
\end{center}

\section{Introduction}
\subsection{Background and motivation}
As an interesting extension of Galton-Watson process, Branching Process in Varying Environments(BPVE hereafter) exhibits some new phenomenons different from Galton-Watson process due to the emerging of the varying environments and has been widely studied in the literature. It is shown that the process may ``fall asleep" at some positive state \cite{lin} and the population may  diverge at different exponential rates \cite{ms}. Moreover, the tail probabilities of the surviving time of the process may show some strange asymptotics \cite{fuj, wy}. For other aspects of the study of BPVEs, we refer the reader to \cite{bcn, cw, dhkp, j, jon, ker,kv} and references therein.

In this paper, we study a BPVE with immigration, namely,  $\{Z_n\}_{n\ge 0}.$ We limit ourselves to BPVEs with geometrical offspring distributions. Our concern is the times the population of the process attaining a fixed level. For $a\in\mathbb Z_+=\{0,1,2,...\},$ let $C=\{n\ge0:Z_n=a\}$ be the collection of times at which the population size of the process attains level $a.$ We aim to give a criterion to determine whether  the set $C$ is finite or not and study further the limit distribution of $|C\cap [1,n]|/\log n$ for the homogeneous Galton-Watson process with immigration, where and in what follows, for a set $A,$ we denote by $|A|$ its cardinality.

 We remark that when $a=0,$ elements in $C$ are usually called regenerating times. Recently in \cite{swyl}, a criterion is given to determine the finiteness of the number of regenerating times for one-type and two-type BPVEs with immigration. Regeneration is also an interesting problem for branching processes in random environments. We refer the reader to \cite{kks} for the single type case and \cite{key,roi,wact} for the multitype case.

 We also remark that there is usually one to one correspondence between questions of branching processes and random walks. So, the results and methods used in this paper are parallel to some extent with those of \cite{w23}, which studies  local time of the nearest-neighbour spatially inhomogeneous random walks on $\mathbb Z_+.$

 Before proceeding further, let us stop to introduce some notations. When  writing $A(n)\sim B(n)$, we mean  $\lim_{n\rto}A(n)/B(n)=1.$ The notation $A(n)=O(B(n))$ means there is a constant $c>0$ such that $|A(n)|<cB(n)$ for all $n$ large enough.  Moreover, for a function $g(x)$ we denote by $g^{(k)}(x)$ its $k$th derivative. We also adopt the convention that empty product equals 1 and empty sum equals $0$. Finally, unless otherwise specified, $0<c<\infty$ is assumed to be a  constant whose value may change from line to line.

\subsection{Models and main results}
Now we introduce precisely the model. For $k\ge1,$ suppose $0<p_k\le 1/2,q_k>0$  are numbers such that  $p_k+q_k=1$ and
 $$f_k(s)=\frac{p_k}{1-q_ks}, s\in [0,1].$$
Let $Z_n,n\ge0$ be a Markov chain such that $Z_0=0$ and
\begin{align}
  E\z( s^{Z_n}\big|Z_0,...,Z_{n-1}\y)=\z[f_{n}(s)\y]^{1+Z_{n-1}}, n\ge1.\label{dz}
\end{align}
Clearly $\{Z_n\}_{n\ge0}$  forms a  BPVE  with exactly one immigrant in each generation.

For $k\ge1,$ let $m_k=f_k'(1)=\frac{q_k}{p_k}.$ Then we have $m_k\ge1, k\ge1.$
For  $n\ge k\ge 1,$ set
\begin{align}
  D(k,n):=1+\sum_{j=k}^n m_j\cdots m_n\text{ and write simply }D(n)\equiv D(1,n). \label{dkn}
\end{align}

Fix $a\in \mathbb Z^+.$ Let $C=\{n\ge0: Z_n=a\}.$ The following theorem gives a criterion for the finiteness of the set $C.$

\begin{theorem}\label{m} Suppose that $p_n\rightarrow 1/2$ as $n\rto.$  Let $D(n), n\ge1$ be the ones in \eqref{dkn}. If $$\sum_{n=2}^\infty\frac{1}{D(n)\log n}<\infty,$$ then $|C|<\infty$ almost surely. If there exists some $\delta>0$ such that $D(n)\le \delta n\log n$ for $n$ large enough and $$\sum_{n=2}^\infty\frac{1}{D(n)\log n}=\infty,$$ then $|C|=\infty$ almost surely.
\end{theorem}
\begin{remark}
  From the viewpoint of Markov chain,  Theorem \ref{m} provides indeed a criterion for the transience and recurrence of the chain $\{Z_n\}_{n\ge0}.$
\end{remark}

Next we give a concrete example.  Fix $B\ge 0$ and let $i_0$ be a positive number such that $\frac{B}{4i_0}<\frac{1}{2}.$ For $i\ge1,$
set \begin{align}\label{pi}
  p_i=\z\{\begin{array}{ll}
    {1}/{2}-\frac{B}{4i}, &i >i_0,\\
    {1}/{2}, &i\le i_0.
  \end{array}\y.
\end{align}
\begin{lemma}\label{del}
  Fix $B\ge0$  and for $i\ge1,$ let $p_i$ be the one in \eqref{pi}. Then
  \begin{align} D(n) \sim \left\{\begin{array}{ll}
             cn^{B}, & B>1, \\
             n\log n, & B=1, \\
             (1-B)^{-1}n, & B<1.
           \end{array}
    \right.\no
    \end{align}
\end{lemma}
 With the help of Lemma \ref{del}, from Theorem \ref{m} we get immediately the following corollary which gives a sharp criterion to determine whether $C$ is a finite set or not.
\begin{corollary}\label{nc}
  Fix $B\ge0$  and for $i\ge1,$ let $p_i$ be the one in \eqref{pi}.
  If $B\ge1,$ then $|C|<\infty$ almost surely. Otherwise, if $B<1,$ then $|C|=\infty$  almost surely.
\end{corollary}
%
%

It is natural to ask how many times the process has ever visited $a$ before time $n.$  We get some partial results in this direction.

\begin{proposition}\label{nr}
  Fix $B\ge 0$  and for $i\ge 1,$ let $p_i$ be the one in \eqref{pi}. Then
   \begin{align}\label{ec}
        \lim_{n\rto}\frac{E |C\cap[1,n]|}{A_n}=1\end{align}
        where
        \begin{align*} A_n=\left\{\begin{array}{ll}
             c, & B>1, \\
             \log\log n, & B=1, \\
             (1-B)\log n, & B<1.
           \end{array}\right.\end{align*}
  \end{proposition}
\proof By Lemma \ref{del} and Lemma \ref{dza} below, we get
  \begin{align}
    P(Z_n=a)=\z(1-\frac{1}{D(n)}\y)^a\frac{1}{D(n)}\sim \frac{1}{D(n)} \sim \left\{\begin{array}{ll}
             cn^{-B}, & B>1, \\
             (n\log n)^{-1}, & B=1, \\
             (1-B)n^{-1}, & B<1.
           \end{array}
    \right.\no
  \end{align}
  As a consequence, we have
  \begin{align*}
    E|C\cap [1,n]|=\sum_{i=0}^n P(Z_i=a)\sim \left\{\begin{array}{ll}
             c, & B>1, \\
             \log\log n, & B=1, \\
             (1-B)\log n, & B<1,
           \end{array}\right. \text{ as } n\rto,
  \end{align*}
which implies \eqref{ec}. The proof is completed. \qed
  \begin{remark}
  From Corollary \ref{nc} and Proposition \ref{nr}, we see that when $B=1,$ though  $|C|<\infty$ almost surely, we have $E|C|=\infty.$ On the other hand, one may expect that
$|C\cap[1,n]|/A_n\rightarrow 1$ almost surely. However, it is not the case in general.

For critical Galton-Watson processes with immigration, we have the following theorem.
\end{remark}

\begin{theorem}\label{ce} If $p_i\equiv \frac{1}{2}$ for all $i\ge1,$ then we have
\begin{align}
  \lim_{n\rto}\frac{|C\cap [1,n]|}{\log n}=S \text{ in distribution, }\label{cd}
\end{align}
  where $S$ is an exponentially distributed random variable with $P(S>t)=e^{-t},\ t>0.$
\end{theorem}

\noindent{\bf Outline of the paper.}
In Section \ref{sec2}, we give some auxiliary results. Section \ref{sec3} is devoted to studying the probabilities  and dependence of the events $\{Z_k=a\}$ and  $\{Z_{k+n}=a\}.$ Finally, Theorem \ref{m} and Theorem \ref{ce} are proved in Section \ref{sec4} and Section \ref{sec5} respectively.

\section{Auxiliary results}\label{sec2}
In this section, we give some preliminary results. To begin with, we give some remarks on  $D(n)$ and $D(k,n), n\ge k\ge1$ which are  defined in \eqref{dkn}.
\subsection{Some facts on $D(n)$ and $D(k,n)$}
%

Notice that $m_i\ge1$ for all $i\ge1.$ Thus we have
\begin{align}
  D(n)>n+1\text{ and } D(k,k+n)>n+2.\label{ubd}
\end{align}
Furthermore, we have $D(n+1)=1+m_{n+1}D(n),$ which implies that
 \begin{align}
D(n+1)>D(n)\text{ for all } n\ge1.\label{dlb}
\end{align}
In addition, some direct computation shows that
\begin{align}
 \frac{D(k,n)}{D(n)}=1-\prod_{j=k-1}^n\z(1-\frac{1}{D(j)}\y).\label{dnp}
\end{align}
The following lemma is useful for studying the dependence between the events $\{k\in C\}$ and $\{k+n\in C\}.$
\begin{lemma}\label{edxy} Suppose that $\lim_{j\rto}m_j=1.$ Then
  for each $\sigma>0$ there exist  $N_1>0$ and $M_1>0$ independent of $k$  such that
  \begin{align}
    1< \frac{D(k,k+n)}{D(k+1,k+n)}< 1+\sigma,   k>N_1,  n>M_1. \label{edkn}
  \end{align}

\end{lemma}
\proof Fix $\sigma>0$ and choose $1>\eta>0$ such that $\frac{2\eta}{(1-\eta)}<\sigma.$ Since $\lim_{j\rto}m_j=1,$
there exists a number $N_1>0$ such that
\begin{align}\label{mi}
  1\ge m_j^{-1}>1-\eta\text{ for all } j>N_1.
\end{align}
It is straightforward to show that
\begin{align}\label{kk}
  \frac{D(k,k+n)}{D(k+1,k+n)}=1+\frac{m_km_{k+1}\cdots m_{k+n}}{D(k+1,k+n)}=1+\frac{1}{\sum_{j=k}^{n+k}m_{k}^{-1}\cdots m_j^{-1}}.
\end{align}
Therefore,
by substituting \eqref{mi} into \eqref{kk},
 for all $n\ge 1, k>N_1$, we get
\begin{align}
  1<\frac{D(k,k+n)}{D(k+1,k+n)}< 1+\frac{1}{\sum_{j=1}^{n+1}(1-\eta)^j}
  =1+\frac{\eta}{(1-\eta)(1-(1-\eta)^{n+1})}.\no
\end{align}
Let $M_1$ be a number(independent of $k$) such that $1-(1-\eta)^{n+1}>1/2$ for all $n>M_1.$ Then we get
\begin{align}
  1<\frac{D(k,k+n)}{D(k+1,k+n)}<
  1+\frac{2\eta}{(1-\eta)}<1+\sigma, k>N_1,n>M_1.\no
\end{align}
The lemma is proved. \qed

Next, we list some combinatorial results which play important roles in the proof of Theorem \ref{ce}.
\subsection{Some combinatorial results}

For positive integers $a\ge j\ge i\ge 1$ set
\begin{align}
  &S(a,j)=\z\{(a_1,...,a_j): a_k\in \mathbb Z_+/\{0\}, j\ge k\ge1, \sum_{k=1}^ja_k=a\y\}.\label{saj}
  \end{align}
The cardinality of $S(a,j)$ is computed in the following lemma. We refer the reader to \cite[Lemma 2]{w23} for its proof.
 \begin{lemma}\label{ns} For $a\ge j\ge i\ge1,$ we have
 $|S(a,j)|=\binom{a-1}{j-1}.$
    \end{lemma}
The next lemma is an improvement of \cite[Lemma 3]{w23}, which will be used for computing the moments of $|C\cap [1,n]|.$
\begin{lemma}\label{lemsa}Fix $k\ge 1,l\ge 1.$ We have
\begin{align}\label{lsc}
  \lim_{n\rto}\frac{1}{(\log n)^k}\sum_{\begin{subarray}{c}l\le j_1<...<j_k\le n,\\
  j_{i+1}-j_i> l,i=1,...,k-1\end{subarray}}\frac{1}{j_1(j_2-j_1)\cdots(j_k-j_{k-1})}=1.
\end{align}
  \end{lemma}
  \proof Notice that  we have (\cite[Lemma 3]{w23})
  \begin{align}\label{lsa}
  \lim_{n\rto}\frac{1}{(\log n)^k}\sum_{l\le j_1<...<j_k\le n}\frac{1}{j_1(j_2-j_1)\cdots(j_k-j_{k-1})}=1.
\end{align}
    Applying \eqref{lsa}, for any $ l_0\ge 1,$ and $1\le i\le k-1,$ we have
    \begin{align}
  \lim_{n\rto}&\frac{1}{(\log n)^k}\sum_{\begin{subarray}{c}l\le j_1<...<j_k\le n,\\
  j_{i+1}-j_i=l_0\end{subarray}}\frac{1}{j_1(j_2-j_1)\cdots(j_k-j_{k-1})}\no\\
  &=\frac{1}{l_0}\lim_{n\rto}\frac{(\log(n-l_0))^{k-1}}{(\log n)^{k}}\no\\
  &\quad\quad\quad\quad\times\frac{1}{(\log (n-l_0))^{k-1}}\sum_{l\le j_1<...<j_{k-1}\le n-l_0}\frac{1}{j_1(j_2-j_1)\cdots(j_{k-1}-j_{k-2})}\no\\
  &=0.\no
\end{align}
Consequently,  we can infer that for any $1\le i\le k-1,$
\begin{align}\label{lsb}
  \lim_{n\rto}&\frac{1}{(\log n)^k}\sum_{\begin{subarray}{c}l\le j_1<...<j_k\le n,\\
  j_{i+1}-j_i\le l\end{subarray}}\frac{1}{j_1(j_2-j_1)\cdots(j_k-j_{k-1})}=0.
\end{align}
Putting \eqref{lsa} and \eqref{lsb} together, we get \eqref{lsc}. Lemma \ref{lemsa} is proved. \qed

\section{Probabilities  and dependence of the events $\{Z_k=a\}$ and $\{Z_{k+n}=a\}$ } \label{sec3}
In this section, we study the distribution of the population size and investigate further the dependence between $\{k\in C\}$ and $\{k+n\in C\}.$ Furthermore, for spatially homogeneous process, we estimate aslo the joint probabilities of $\{j_i\in C\},1\le j_1<...<j_m\le n.$

\subsection{Distribution of the population size}
To begin with, we compute firstly the distribution of $Z_n.$
For $n\ge0,$ let
 $F_n(s)=E(s^{Z_n}),s\in [0,1].$ By induction, from \eqref{dz}  we get
\begin{align}\label{fn}
  F_n(s)=\prod_{k=1}^n f_{k,n}(s), s\in [0,1], n\ge0,
\end{align}
where
$f_{k,n}(s)=f_k(f_{k+1}(\cdots(f_n(s)) )), n\ge k\ge 1.$
Recall that $m_k=f_k'(1)=\frac{q_k}{p_k},k\ge1.$ Then it is easily seen that
$f_k(s)=1-\frac{m_k(1-s)}{1+m_k(1-s)}, k\ge1.$ As a consequence, it follows by induction that
\begin{align}
  f_{k,n}(s)&=1-\frac{m_k\cdots m_n(1-s)}{1+\sum_{j=k}^n m_j\cdots m_n(1-s)}
  =\frac{1+\sum_{j=k+1}^n m_j\cdots m_n(1-s)}{1+\sum_{j=k}^n m_j\cdots m_n(1-s)}\no\\
  &=\frac{D(k+1,n)+(D(k+1,n)-1)s}{D(k,n)+(D(k,n)-1)s}.\label{fkn}
\end{align}
If we plug \eqref{fkn} back into \eqref{fn}, we have
\begin{align*}
  F_n(s)=\frac{1}{1+\sum_{j=1}^n m_j\cdots m_n(1-s)}.
\end{align*}
By some easy manipulation, we get
\begin{align}
  F^{(a)}_n(s)=\frac{a!\z(\sum_{j=1}^n m_j\cdots m_n\y)^a}{\z(1+\sum_{j=1}^n m_j\cdots m_n(1-s)\y)^{a+1}},s\in[0,1].\no
\end{align}
As a consequence, we obtain
\begin{align}
  P(Z_n=a)&=\frac{F^{(a)}_n(0)}{a!}=\frac{\z(\sum_{j=1}^n m_j\cdots m_n\y)^a}{\z(1+\sum_{j=1}^n m_j\cdots m_n\y)^{a+1}}\no\\
  &=\frac{(D(n)-1)^a}{D(n)^{a+1}}=\z(1-\frac{1}{D(n)}\y)^a\frac{1}{D(n)}.\no
\end{align}
To summarize,  we have the following lemma.
\begin{lemma}\label{dza}
  For $n\ge1,$ we have \begin{align}\label{ha}
  &P(Z_n=a)=\frac{(D(n)-1)^a}{D(n)^{a+1}}=\z(1-\frac{1}{D(n)}\y)^a\frac{1}{D(n)}.
\end{align}
\end{lemma}

\subsection{Dependence between $\{Z_k=a\}$ and $\{Z_{k+n}=a\}$}
An important step to prove Theorem \ref{m} is the characterize the dependence between the events $\{Z_k=a\}$ and $\{Z_{k+n}=a\}$ which is done in the following proposition.
\begin{proposition}\label{dpkn}  For each $\ve>0,$ there exist  $N>0$ and $M>0$ independent of $k$ such that
\begin{align}
 1-\ve\le \frac{D(k+1,k+n)}{D(k+n)}\frac{P(Z_{k+n}=a,Z_{k}=a)}{P(Z_{k+n}=a)P(Z_{k}=a)}
  \le 1+\ve, k>N,n>M.\label{je}
\end{align}

\end{proposition}
\proof Fix $\ve>0$ and choose $\sigma>0$ such that $1-\ve<(1-\sigma)^2<(1+\sigma)^2<1+\ve.$  In order to prove \eqref{je},
for $n,k\ge1$ we define
\begin{align}
  G_{nk}(s|a)=E(s^{Z_{n+k}}|Z_k=a),s\in[0,1].\no
\end{align}
Clearly, owing to \eqref{fkn}, we have
\begin{align}
  G_{nk}(s|a)&=\z(f_{k+1,k+n}(s)\y)^a\prod_{j=k+1}^{k+n}f_{j,k+n}(s)\no\\
  &=\z(\frac{D(k+2,k+n)-(D(k+2,k+n)-1)s}{D(k+1,k+n)-(D(k+1,k+n)-1)s}\y)^a\no\\
  &\quad\quad\quad\quad\times\frac{1}{D(k+1,k+n)-(D(k+1,k+n)-1)s}\no\\
  &=:B(s)^aA(s)=:g(s)A(s),\label{gnksa}
\end{align}
where we temporarily write $A(s)=\frac{1}{D(k+1,k+n)-(D(k+1,k+n)-1)s}$  and $g(s)=B(s)^a$ with
$B(s)=\frac{D(k+2,k+n)-(D(k+2,k+n)-1)s}{D(k+1,k+n)-(D(k+1,k+n)-1)s}.$
Therefore, we get
\begin{align}
  P(Z_{k+n}&=a|Z_{k}=a)=\frac{G_{nk}^{(a)}(0|a)}{a!}
  =\frac{1}{a!}\sum_{j=0}^a \binom{a}{j}A^{(a-j)}(0)g^{(j)}(0)\no\\
&=\frac{1}{a!}A^{(a)}(0)g(0)+\frac{1}{a!}\sum_{j=1}^{(a)}A^{(a-j)}(0)g^{(j)}(0)\no\\
 &=\frac{(D(k+1,k+n)-1)^a}{D(k+1,k+n)^{a+1}}\z(\frac{D(k+2,k+n)}{D(k+1,k+n)}\y)^a\no\\
 &\quad\quad\quad\quad\quad\quad\quad\quad\quad\quad+\sum_{j=1}^{a}\frac{1}{j!}\frac{(D(k+1,k+n)-1)^{a-j} }{D(k+1,k+n)^{a-j+1}}g^{(j)}(0)\no\\
 &=\frac{(D(k+1,k+n)-1)^a}{D(k+1,k+n)^{a+1}}\no\\
 &\quad\quad\quad\quad\times\z(\z(\frac{D(k+2,k+n)}{D(k+1,k+n)}\y)^a+
 \sum_{j=1}^{a}\frac{1}{j!}\frac{D(k+1,k+n)^j }{(D(k+1,k+n)-1)^j }g^{(j)}(0)\y).\label{knu}
\end{align}
We claim that
\begin{align}\label{gu}
  g^{(j)}(0)\rightarrow 0 \text{ uniformly in } k\text{ as } n\rto.
\end{align}
In fact, note that $B(0)=\frac{D(k+2,k+n)}{D(k+1,k+n)}<1$ and   for $1\le j\le a$ taking \eqref{ubd} into consideration, we have
\begin{align}
  B^{(j)}(0)&=j!(D(k+1,k+n)-D(k+2,k+n))
  \frac{(D(k+1,k+n)-1)^{j-1}}{(D(k+1,k+n))^{j+1}}\no\\
  &\le j!\frac{m_{k+1}\cdots m_{k+n}}{D(k+1,k+n)^2}
  \le j!\frac{1}{D(k+1,k+n)}\rightarrow 0\no
\end{align}
uniformly in $k\ge 0$ as $n\rto.$  Because $g^{(j)}(0)$ is the sum of the product of the form $i\binom{j}{l}B(0)^{a-l}\prod_{h=1}^{j}\z(B^{(h)}(0)\y)^{b_h}$
where $1\le i,l,h\le j,$  $1\le b_h\le j,$ we infer  that \eqref{gu} is true.

Putting \eqref{ubd}, \eqref{edkn} and \eqref{gu} together, we come to a conclusion that there exist $N_2>0$ and $M_2>0$ independent of $k$ such that
\begin{align}
  1-\sigma<\z(\frac{D(k+2,k+n)}{D(k+1,k+n)}\y)^a+
 \sum_{j=1}^{a}\frac{1}{j!}\frac{D(k+1,k+n)^j }{(D(k+1,k+n)-1)^j }g^{(j)}(0)<1+\sigma\no
\end{align}
for all $k>N_2,n>M_2.$ Moreover,
it follows from \eqref{ubd} that for some $M_3>0$ independent of $k,$ we have
\begin{align}
1<\z(\frac{D(k+n)}{D(k+n)-1}\y)^a<1+\ve,\  1-\ve< \z(\frac{D(k+1,k+n)-1}{D(k+1,k+n)}\y)^a<1, n>M_3.\no
\end{align}
Notice that by \eqref{ha}, we have
$$
\frac{D(k+1,k+n)}{D(k+n)}\frac{P(Z_{k+n}=a,Z_{k}=a)}{P(Z_{k+n}=a)P(Z_{k}=a)}=\frac{P(Z_{k+n}=a|Z_{k}=a)(D(k+n))^{a}}{(D(k+1,k+n))^{-1}(D(k+n)-1)^{a}}.
  $$
Consequently, setting $N=N_2$ and  $M=M_2\vee M_3,$ we infer from \eqref{ubd}, \eqref{edkn} and \eqref{knu} that
\begin{align}
(1-\sigma)^2 <\frac{D(k+1,k+n)}{D(k+n)}\frac{P(Z_{k+n}=a,Z_{k}=a)}{P(Z_{k+n}=a)P(Z_{k}=a)}\le (1+\sigma)^2,k>N, n>M.\no
\end{align}
The proposition is proved. \qed

\subsection{Homogeneous setting}
In this subsection, we investigate the homogeneous setting, that is, the case $p_n\equiv1/2$ for all $i\ge1.$   As a byproduct, we also get a combinatorial formula which may have its own interest. The following lemma is the main result of this subsection.
 \begin{lemma}\label{hs} For $a\in \mathbb Z_+,$ we have
   \begin{align}
  \sum_{j=0}^a\frac{(a+j)!}{j!j!(a-j)!}(-1)^{a-j}=1.\label{l2}
\end{align}
Furthermore, if $p_n\equiv1/2$ for all $i\ge1,$ we have for $n\ge 1, k\ge 0,$
\begin{align}\label{fa}
  P(Z_{n+k}=a|Z_k=a)
  =\sum_{j=0}^a\frac{(a+j)!}{j!j!(a-j)!}\z(\frac{n}{n+1}\y)^j\z(\frac{n}{n+1}\y)^j\z(\frac{1-n}{n+1}\y)^{(a-j)}\frac{1}{n+1},
\end{align}
and
\begin{align}
  \lim_{n\rto}nP(Z_{k+n}=a|Z_k=a)=1,k\ge 0.\label{pah}
\end{align}
 \end{lemma}
 \proof
  Notice that when $p_n\equiv1/2$ for all $i\ge1,$ we have
  \begin{align}
    D(n)=n+1, D(k,n)=n-k+2, n\ge k\ge1.\label{dh}
  \end{align}
  Again we write
  \begin{align}
  G_{nk}(s|a)=E(s^{Z_{n+k}}|Z_k=a),s\in[0,1].\no
\end{align}
  Substituting \eqref{dh} into \eqref{gnksa}, we get
  \begin{align}\no
  G_{nk}(s|a)=\z(\frac{n-(n-1)s}{n+1-ns}\y)^a\frac{1}{n+1-ns}.
\end{align}
Furthermore, taking \eqref{knu} and \eqref{gu} together, we deduce that
\begin{align}\label{lgn}
  \lim_{n\rto}nP(Z_{k+n}=a|Z_k=a)=\lim_{n\rto}\frac{nG_{nk}^{(a)}(0|a)}{a!}=1,k\ge 0
\end{align} which proves \eqref{pah}.
But from \eqref{knu}, we can not obtain an explicit formula for $P(Z_{k+n}=a|Z_k=a).$
Therefore we  write
$$G_{nk}(s|a)=C(s)D(s)$$ where
\begin{align*}
  C(s)=(n-(n-1)s)^a, D(s)=(n+1-ns)^{-(a+1)}.
\end{align*}
By some routine computation, for $a\ge j\ge0,$ we get
\begin{align*}
  C^{(j)}(0)&=a(a-1)\cdots (a-j+1)n^{a-j}(1-n)^{j},\\
  D^{(j)}(0)&=(a+1)(a+2)\cdots (a+j)(n+1)^{-(a+j+1)}n^j.
\end{align*}
As a result, we get
\begin{align}
  G_{nk}^{(a)}(0|a)&=\sum_{j=0}^a\binom{a}{j}C^{(a-j)}(0)D^{(j)}(0)\no\\
  &=\sum_{j=0}^a\binom{a}{j}\frac{(a+j)!}{j!}n^{j}(1-n)^{a-j}(n+1)^{-(a+j+1)}n^j\no\\
  &=\sum_{j=0}^a\frac{a!}{j!(a-j)!}\frac{(a+j)!}{j!}\z(\frac{n}{n+1}\y)^j\z(\frac{n}{n+1}\y)^j\z(\frac{1-n}{n+1}\y)^{(a-j)}\frac{1}{n+1},\no
\end{align} which implies \eqref{fa} and \eqref{pah}.
On the other hand, we infer from \eqref{pah} that
\begin{align}
  \lim_{n\rto}\frac{n}{a!}G_{nk}^{(a)}(0|a)=\sum_{j=0}^a\frac{(a+j)!}{j!j!(a-j)!}(-1)^{a-j}.\label{l3}
\end{align}
Taking \eqref{lgn} and \eqref{l3} together, we obtain \eqref{l2}. The proof is completed. \qed

To close this subsection, we give the following lemma which will be used to prove Theorem \ref{ce}.

\begin{lemma}\label{hej}
  Fix $m\ge 2.$ For $\ve>0,$ there exists a number $n_0>0$ such that
  \begin{align*}
    (1-\ve)^m< \frac{P(Z_{j_1}=a,...,Z_{j_m}=a)}{j_1^{-1}\prod_{i=1}^{m-1}(j_{i+1}-j_i)^{-1}}<(1+\ve)^m,
  \end{align*}
  for  $j_1>n_0, j_{i+1}-j_i>n_0,i=1,...,m-1.$
  \end{lemma}
\proof Fix $\ve>0.$ By \eqref{pah}, there exists a number $n_1$ independent of $k$ such that $1-\ve<nP(Z_{n+k}=a|Z_k=a)<1+\ve$ for all $n\ge n_1,k\ge1.$ Therefore we  have \begin{align}1-\ve<P(Z_{j_{i+1}}=a|Z_{j_i}=a)(j_{i+1}-j_i)<1+\ve,\ j_{i+1}-j_i>n_1.\label{pc}\end{align} Moreover from \eqref{ha} and \eqref{dh}, we see that for $n\ge1,$ $P(Z_{n}=a)=\z(1-\frac{1}{n+1}\y)^a\frac{1}{n+1}.$ Thus, there exists a number $n_2$ such that \begin{align}
  1-\ve<j_1 P(Z_{j_1}=a)<1+\ve,\ j_1>n_2.\label{pj1}
\end{align}
Note that
\begin{align*}
P(Z_{j_1}=a,...,Z_{j_m}=a)=P(Z_{j_1}=a)\prod_{i=1}^{m-1}P(Z_{j_{i+1}}=a|Z_{j_i}=a).
\end{align*}
Then letting $n_0=n_1\vee n_2$  we get
  \begin{align*}
    (1-\ve)^m< \frac{P(Z_{j_1}=a,...,Z_{j_m}=a)}{j_1^{-1}\prod_{i=1}^{m-1}(j_{i+1}-j_i)^{-1}}<(1+\ve)^m,
  \end{align*}
if $j_1>n_0$ and $j_{i+1}-j_i>n_0,i=1,...,m-1.$ The lemma is proved. \qed


\section{Criteria of the finiteness of $C$-Proof of Theorem \ref{m}}\label{sec4}
 We prove Theorem \ref{m} by an approach similar to the one used in \cite{cfrb}. We divide the proof into two parts, that is, the convergent one and the divergent one.

 \subsection{Convergent part}
 \proof  To prove the convergent part, we follow the ideas of \cite{jlp}. For $j<i,$ denote  $C_{j,i}=\{(2^j, 2^i]: x\in C\}$ and set $A_{j,i}=|C_{j,i}|.$ On the event $\{A_{m,m+1}>0\},$  let $l_m=\max\{k:k\in C_{m,m+1}\}.$
Then for $m\ge1,$ we have
\begin{align}
  &\sum_{j=2^{m-1}+1}^{2^{m+1}}P(j\in C)=E(A_{m-1,m+1})\no\\
    &\quad\quad\ge \sum_{n=2^m+1}^{2^{m+1}}E(A_{m-1,m+1}, A_{m,m+1}>0,l_m=n)\no\\
    &\quad\quad=\sum_{n=2^m+1}^{2^{m+1}}P(A_{m,m+1}>0,l_m=n)E(A_{m-1,m+1}| A_{m,m+1}>0,l_m=n)\no\\
    &\quad\quad=\sum_{n=2^m+1}^{2^{m+1}}P(A_{m,m+1}>0,l_m=n)\sum_{i=2^{m-1}+1}^nP(i\in C| A_{m,m+1}>0,l_m=n)\no\\
    &\quad\quad\ge P(A_{m,m+1}>0)\min_{2^m<n\le 2^{m+1}}\sum_{i=2^{m-1}+1}^nP(i\in C| A_{m,m+1}>0,l_m=n)\no\\
    &\quad\quad=:a_mb_m. \label{smpjc}
\end{align}
Fix $2^{m}+1\le n\le 2^{m+1}$ and $2^{m-1}+1\le i\le n.$ Using Markov property, we get
\begin{align}
  P(&i\in C| A_{m,m+1}>0,l_m=n)\no\\
  &=\frac{P\z(Z_i=a, Z_n=a, Z_t\ne a, n+1\le t\le 2^{m+1}\y)}{P\z(Z_n=a, Z_t\ne a, n+1\le t\le 2^{m+1}\y)}\no\\
  &=\frac{P(Z_i=a,Z_n=a)}{P(Z_n=a)}\frac{P\z( Z_t\ne a, n+1\le t\le 2^{m+1}|Z_i=a, Z_n=a\y)}{P\z(Z_t\ne a, n+1\le t\le 2^{m+1}|Z_n=a\y)}\no\\
  &= \frac{P(Z_i=a,Z_n=a)}{P(Z_n=a)}=\frac{P(Z_i=a)P(Z_n=a|Z_i=a)}{P(Z_n=a)}.\label{nl}
\end{align}
Since $\lim_{k\rto}m_k=1,$ $\sup_{k\ge1}m_k<B$ for some constant $B>0.$ Thus,
taking \eqref{kk} into account, from \eqref{knu}, we see that
\begin{align}
  P(Z_{n}=a|Z_{i}=a)&\ge\z(1-\frac{1}{D(i+1,n)}\y)^a\z(\frac{D(i+2,n)}{D(i+1,n)}\y)^a\frac{1}{D(i+1,n)}\no\\
&\ge \z(1-\frac{1}{1+n-i}\y)^a\z(\frac{1}{1+m_{i+1}}\y)^a\frac{1}{D(i+1,n)}\no\\
&\ge  2^{-a}(1+B)^{-a}\frac{1}{D(i+1,n)},\ i<n.\label{cl}
\end{align}
Using the fact $1/D(i)<1/2$ for all $i,$ we can infer from \eqref{ha}, \eqref{nl} and \eqref{cl} that
\begin{align}
   P(i\in C| A_{m,m+1}>0,l_m=n)&\ge 2^{-a}(1+B)^{-a} \frac{\z(1-\frac{1}{D(i)}\y)^a\frac{1}{D(i)}}{\z(1-\frac{1}{D(n)}\y)^a\frac{1}{D(n)}}
\frac{1}{D(i+1,n)}\no\\
&\ge 4^{-a}(1+B)^{-a}\frac{D(n)}{D(i)D(i+1,n)}, i<n.\label{pnil}
\end{align}
  As $m_j\ge 1$ for all $j\ge1,$ we have
  \begin{align}
    &\frac{D(n)}{D(i)D(i+1,n)}=\frac{\sum_{j=1}^{n+1}m_j\cdots m_n}{\sum_{j=1}^{i+1}m_j\cdots m_i\sum_{j=i+1}^{n+1}m_j\cdots m_n}\no\\
    &\quad\quad=\frac{\sum_{j=1}^{n+1}m_1^{-1}\cdots m_{j-1}^{-1}}{\sum_{j=1}^{i+1}m_1^{-1}\cdots m_{j-1}^{-1}\sum_{j=i+1}^{n+1}m_{i+1}^{-1}\cdots m_{j-1}^{-1}}\no\\
    &\quad\quad\ge \frac{1}{\sum_{j=i+1}^{n+1}m_{i+1}^{-1}\cdots m_{j-1}^{-1}}\ge\frac{1}{n-i+1}.\label{pd}   \end{align}
Taking \eqref{pnil} and \eqref{pd} together, we deduce that
\begin{align}
  b_m&=\min_{2^m<n\le 2^{m+1}}\sum_{i=2^{m-1}+1}^nP(i\in C| A_{m,m+1}>0,l_m=n)\no\\
  &\ge  4^{-a}(1+B)^{-a}\min_{2^m<n\le 2^{m+1}}\sum_{i=2^{m-1}+1}^n \frac{1}{n-i+1}\no\\
  &= 4^{-a}(1+B)^{-a}\min_{2^m<n\le 2^{m+1}}\sum_{j=1}^{n-2^{m-1}} \frac{1}{j}\no\\
  &= 4^{-a}(1+B)^{-a}\sum_{j=1}^{2^{m-1}+1} \frac{1}{j}\ge 4^{-a}(1+B)^{-a}\int_{1}^{2^{m-1}+2}\frac{1}{x}dx\no\\
  & \ge 4^{-a}(1+B)^{-a} (m-1)\log 2. \label{ebm}
\end{align}
Substituting \eqref{ebm} into \eqref{smpjc} and using \eqref{ha}, we see that
\begin{align}
  \sum_{m=2}^{\infty}&P(A_{m,m+1}>0)\le  \sum_{m=2}^{\infty}\frac{1}{b_m}\sum_{j=2^{m-1}+1}^{2^{m+1}}P(j\in C)\no\\
  &\le \frac{4^a(1+B)^a}{\log 2}\sum_{m=2}^{\infty} \frac{1}{m-1}\sum_{j=2^{m-1}+1}^{2^{m+1}}\z(1-\frac{1}{D(j)}\y)^a\frac{1}{D(j)}\no\\
  &\le \frac{4^a(1+B)^a}{\log 2}\sum_{m=2}^{\infty} \sum_{j=2^{m-1}+1}^{2^{m+1}}\frac{1}{D(j)\log j}
  \le \frac{4^a(1+B)^a}{\log 2}\sum_{n=2}^\infty \frac{1}{D(n)\log n}.\no
\end{align}
Therefore, if $\sum_{n=2}^\infty \frac{1}{D(n)\log n}<\infty,$ then it follows by the Borel-Cantelli lemma that with probability 1, at most finitely many of the events $\{A_{m,m+1}>0\},m\ge1$ occur. Consequently, $|C|<\infty$ almost surely. The convergent part of Theorem \ref{m} is proved. \qed

\subsection{Divergent part}
\proof Suppose there exists some $\delta>0$ such that $D(n)\le \delta n\log n$ for $n$ large enough and $\sum_{n=2}^\infty\frac{1}{D(n)\log n}=\infty.$
For $j\ge1,$ let $n_j=[j\log j]$ be the integral part of $j\log j$ and set $A_j=\{n_j\in C\}.$
From \eqref{ha}, we get
\begin{align}
\sum_{j=2}^{\infty}P(A_j)=\sum_{j=2}^{\infty}\z(1-\frac{1}{D(n_j)}\y)^a\frac{1}{D(n_j)}\ge 2^{-a}\sum_{j=2}^{\infty}\frac{1}{D([j\log j])}. \label{ad}
\end{align}
Since $m_j\ge 1$ for all $j\ge1,$ from \eqref{dlb} we see that $D(n)$ is increasing in $n.$
Thus it follows from \cite[Lemma 2.2]{cfrb} that
$\sum_{j=2}^{\infty}\frac{1}{D([j\log j])}$ and $\sum_{j=2}^{\infty}\frac{1}{D(j)\log j}$  converge or diverge simultaneously.  Therefore,  from \eqref{ad} we see that
\begin{align}
  \sum_{j=2}^{\infty}P(A_j)=\infty.\label{spa}
\end{align}

Next fix $\varepsilon >0.$
Obviously,
 for $l>k,$
$n_l-n_k\ge l\log l-k\log k\ge \log k.$ Let $M,N$ be those in \eqref{je}. Then, there exists a number $N_3>0$ such that for all $k\ge N_3,$
$n_k\ge N, n_l-n_k>M.$
Therefore, applying   \eqref{je} and taking \eqref{dnp} into account, for $l>k>N_3,$ we have
\begin{align}\label{e3}
P(A_kA_l)
&\leq (1+\varepsilon)P(A_k)P(A_l)\z(\frac{D(n_k,n_l)}{D(n_k)}\y)^{-1}\nonumber\\
&=(1+\varepsilon)\z(1-\prod_{i=n_k-1}^{n_l}\z(1-\frac{1}{D(i)}\y)\y)^{-1}P(A_k)P(A_l)\no\\
&\le(1+\varepsilon)\z(1-\exp\z\{-\sum_{i=n_k-1}^{n_l}\frac{1}{D(i)}\y\}\y)^{-1}P(A_k)P(A_l).
\end{align}
Define
\begin{align}
\ell_1=\min\z\{l\geq k:\sum_{i=n_k-1}^{n_l}\frac{1}{D(i)}\geq \log \frac{1+\varepsilon}{\varepsilon}\y\}.\no
\end{align}
Then for $l\geq \ell_1$ we have
$\z(1-\exp\z\{-\sum_{i=n_k-1}^{n_l}\frac{1}{D(i)}\y\}\y)^{-1}\leq 1+\varepsilon.$
Therefore, from (\ref{e3}) we get
\begin{align}
P(A_kA_l)\leq (1+\varepsilon)^2P(A_k)P(A_l), \text{  for } l\geq \ell_1.\label{e6}
\end{align}

Next we consider $k<l<\ell_1.$ Note that for $0\leq u\leq \log \frac{1+\varepsilon}{\varepsilon}$, we have $1-e^{-u}\geq cu$. Since $D(k)$ is increasing in $k\ge1,$  thus on account of \eqref{ha}, and (\ref{e3}), we have
\begin{align}
P(A_kA_l)&\le(1+\varepsilon)\z(1-\exp\z\{-\sum_{i=n_k-1}^{n_l}\frac{1}{D(i)}\y\}\y)^{-1}P(A_k)P(A_l)\nonumber\\
&\leq\frac{cP(A_k)P(A_l)}{\sum_{i=n_k-1}^{n_l}\frac{1}{D(i)}}
\leq \frac{cP(A_k)P(A_l)D(n_l)}{n_l-n_k+2}\le \frac{cP(A_k)}{l\log l-k\log k}.\no
\end{align}
As a consequence, we get
\begin{align}
\sum_{l=k+1}^{\ell_1-1}&P(A_kA_l)\leq cP(A_k)\sum_{l=k+1}^{\ell_1-1}\frac{1}{l\log l-k\log k}\nonumber\\
&\leq cP(A_k)\frac{1}{\log k}\sum_{l=k+1}^{\ell_1-1}\frac{1}{l-k}\leq cP(A_k)\frac{\log \ell_1}{\log k}\le cP(A_k),\label{e8}
\end{align}where for the last inequality, we use the fact
$
\frac{\log \ell_1}{\log k}\leq \gamma$
with some constant $\gamma$ depending only on $\varepsilon,$ see \cite[p.635, display (4.2)]{cfrb}.

Taking (\ref{e6}) and (\ref{e8}) together, we have
$$\sum_{k=N_3}^{n}\sum_{l=k+1}^{n}P(A_kA_l)\leq\sum_{k=N_3}^{n}\sum_{l=k+1}^{n}(1+\varepsilon)^2P(A_k)P(A_l)+c\sum_{k=N_3}^{n}P(A_k).$$
Writing $H(\varepsilon)=(1+\varepsilon)^2$, thanks to (\ref{spa}), we have
\begin{align*}
\alpha_{H}&:=\liminf_{n\rightarrow\infty}\frac{\sum_{k=N_3}^{n}\sum_{l=k+1}^{n}P(A_kA_l)-\sum_{k=N_3}^{n}\sum_{l=k+1}^{n}HP(A_k)P(A_l)}{[\sum_{k=N_3}^{n}P(A_k)]^2}\nonumber\\
&\leq\lim_{n\rightarrow\infty}\frac{c}{\sum_{k=N_3}^{n}P(A_k)}=0.
\end{align*}
Applying  a generalized version of Borel-Cantelli lemma \cite[p.235]{pe04}, we have
$$P(A_k,k\geq N_3 \text{ occur infinitely often})\geq\frac{1}{H+2\alpha_{H}}\ge\frac{1}{(1+\varepsilon)^2}.$$
Thus, we obtain 
\begin{align}
  P(A_k,k\geq 1 \text{ occur infinitely often})\ge P(A_k,k\geq N_3 \text{ occur infinitely often}) \ge\frac{1}{(1+\varepsilon)^2}.\no
\end{align}
As $\varepsilon$ is arbitrary, letting $\varepsilon\rightarrow0,$ we conclude that
$$P(A_k,k\geq 1 \text{ occur infinitely often})=1.$$ This completes the proof of the divergent part of Theorem \ref{m}. \qed

\section{The cardinality of $C\cap[1,n]$}\label{sec5}
In this section, we study  asymptotics of the cardinality $C\cap[1,n],$ giving the proofs of Theorem \ref{ce}. To begin with, we give the proof of Lemma \ref{del}.

\subsection{Proof of Lemma \ref{del}}
\proof
Fix $B\ge0$ and for $i\ge 1,$ let $p_i$ be the one in \eqref{pi}. Then we have $$m_i=\frac{q_i}{p_i}=\frac{\frac{1}{2}+\frac{B}{4i}}{\frac{1}{2}-\frac{B}{4i}}=1+\frac{B}{i}+O\z(\frac{1}{i^2}\y)$$ as $i\rto.$ Consequently, for some $c_0>0,$ we have
\begin{align}\label{ms}
  m_1\cdots m_n\sim c_0^{-1} n^B, \text{ as }n\rto.
\end{align}
With \eqref{ms} in hand, by some easy computation, we can show that
\begin{align}
 1+ \sum_{i=1}^{n}m_1^{-1}\cdots m_i^{-1}\sim \left\{\begin{array}{ll}
    c_1,&\text{ if }B>1,\\
    c_0\log n,& \text{ if }B=1,\\
    c_0(1-B)^{-1}n^{1-B},&\text{ if }B<1,
  \end{array}\right. \text{ as }n\rto,\no
\end{align}
where $c_1>0$ is a proper constant.
As a result, we have
\begin{align}\no
   D(n)=\frac{1+\sum_{i=1}^{n}m_1^{-1}\cdots m_{i}^{-1}}{m_1^{-1}\cdots m_n^{-1}}\sim \left\{\begin{array}{ll}
    c_1/c_0 n^{B},&\text{ if }B>1,\\
    n\log n,& \text{ if }B=1,\\
    (1-B)^{-1}n,&\text{ if }B<1,
  \end{array}\right.
\end{align}
which  completes the proof of Lemma \ref{del}. \qed

\subsection{Proof of Theorem \ref{ce}}

\proof To begin with we show that the moments of  $|C\cap [1,n]|$ converges to those of exponential distribution, that is
 \begin{align}
   \lim_{n\rto}\frac{ E(|C\cap [1,n]|^k)}{(\log n)^k}=k!,\ k\ge1.\label{cme}
 \end{align}
 In fact, setting $\eta_k=\z\{\begin{array}{cc}
              1, & k\in C \\
              0 &  k\notin C
            \end{array}
\y.$ for $k\ge1,$  then we have $|C\cap [1,n]|=\sum_{k=1}^n \eta_k.$
Therefore, with $S(k,m)$ the one defined in \eqref{saj}, we have 
\begin{align}
  E|C\cap [1,n]|^k&=E\z(\z(\sum_{j=1}^n\eta_j\y)^k\y)=\sum_{1\le j_1,j_2,...,j_k\le n}E(\eta_{j_1}\eta_{j_2}\cdots\eta_{j_k})\no\\
  &=\sum_{m=1}^k\sum_{\begin{subarray}{c}l_1+\dots+l_m=k,\\
  l_i\ge1,i=1,...,m
  \end{subarray}}m!\sum_{0\le j_1<...<j_m\le n}E(\eta_{j_1}^{l_1}\cdots\eta_{j_m}^{l_m})\no\\
  &=\sum_{m=1}^k\sum_{(l_1,...,l_m)\in S(k,m)}m!\sum_{0\le j_1<...<j_m\le n}E(\eta_{j_1}^{l_1}\cdots\eta_{j_m}^{l_m}).\no
\end{align}
Since the values of $\eta_j,j=1,...,n$ are either  $0$ or $1,$ then on account of Lemma \ref{ns} we have
\begin{align}
  E|C\cap [1,n]|^k&=\sum_{m=1}^k\binom{k-1}{m-1}m!\sum_{1\le j_1<...<j_m\le n}P(Z_{j_1}=a,\cdots ,Z_{j_m}=a)\no\\
  &=:\sum_{m=1}^k\binom{k-1}{m-1}m!G(m,n).\label{eck}
  \end{align}
  Next we show that
  \begin{align}
    \lim_{n\rto}\frac{G(m,n)}{(\log n)^m}=1,1\le m\le k.\label{lgmn}
  \end{align}
    To this end, we fix $\ve>0,$  let $n_0$ be the one in Lemma \ref{hej} and   denote temporarily \begin{align*}
    A&=\{(j_1,...,j_m)\mid 1\le j_1<....<j_m\le n, j_{i+1}-j_i>n_0, \forall 0\le i\le m-1\},\\
    B&=\{(j_1,...,j_m)\mid 1\le j_1<....<j_m\le n, (j_1,...,j_m)\notin A\},
  \end{align*} where and in the remainder of this proof we set $j_0=0.$
  Then we can write
    \begin{align}
    G(m,n)&=\sum_{(j_1,...,j_m)\in A}+\sum_{(j_1,...,j_m)\in B}
P(Z_{j_1}=a,\cdots ,Z_{j_m}=a)\no\\
&=:G_A(m,n)+G_{B}(m,n).\label{gab}
  \end{align}
  We claim that
  \begin{align}
    & (1+\ve)^m<\varliminf_{n\rto}\frac{G_A(m,n)}{(\log n)^m}\le \varlimsup_{n\rto}\frac{G_A(m,n)}{(\log n)^m}\le (1+\ve)^m,\label{gmn1}
    \end{align}
    and
    \begin{align}
    &\lim_{n\rto}\frac{G_B(m,n)}{(\log n)^m}=0.\label{gmn0}
  \end{align}
  To prove \eqref{gmn1}, notice that by Lemma \ref{hej}
  we have
  \begin{align}
   \sum_{(j_1,...,j_m)\in A}&\frac{(1-\ve)^m}{j_1\prod_{i=1}^{m-1}(j_{i+1}-j_i)} \le G_A(m,n)\le \sum_{(j_1,...,j_m)\in A}\frac{(1+\ve)^m}{j_1\prod_{i=1}^{m-1}(j_{i+1}-j_i)}.\label{gau}
  \end{align}
  Thus, dividing \eqref{gau} by $(\log n)^m$ and taking the upper and the lower limits respectively,  in view of Lemma \ref{lemsa}, we get  \eqref{gmn1}.

  Next we prove \eqref{gmn0}. For this purpose, for
  $0\le i_1<...<i_l\le m-1,$ and $k_1,...,k_l\ge1,$  set
  $$B(i_1,k_1,...,i_l,k_l)=\z\{(j_1,...,j_m)\z|\begin{array}{c}1\le j_1<....<j_m\le n, \\
    j_{i_s+1}-j_{i_s}=k_s,0\le s\le l\\
    j_{i+1}-j_i>n_0, i\in\{0,..,m-1\}/\{i_1,...,i_l\}\end{array}\y.\y\}.$$
  Then we can write
  \begin{align*}
    B&=\bigcup_{l=1}^{m-1}\bigcup_{\begin{subarray}{c}
                                    0\le i_1<...<i_l\le m-1\\
                                     1\le k_1,...,k_l\le n_0
                                  \end{subarray}
    }\z\{(j_1,...,j_m)\z|\begin{array}{c}1\le j_1<....<j_m\le n,j_{i_s+1}-j_{i_s}=k_s, \\
0\le s\le l, j_{i+1}-j_i>n_0,\\
     i\in\{0,..,m-1\}/\{i_1,...,i_l\}\end{array}\y.\y\}\no\\
   & =\bigcup_{l=1}^{m-1}\bigcup_{\begin{subarray}{c}
                                    0\le i_1<...<i_l\le m-1\\
                                     1\le k_1,...,k_l\le n_0
                                  \end{subarray}
    }B(i_1,k_1,...,i_l,k_l).
  \end{align*}
  Therefore, we have
  \begin{align}
    G_B(m,n)&=\sum_{(j_1,...,j_m)\in B}
P(Z_{j_1}=a,\cdots ,Z_{j_m}=a)\no\\
&=\sum_{(j_1,...,j_m)\in B}P(Z_{j_1}=a)\prod_{i=1}^{m-1}P(Z_{j_{i+1}}=a|Z_{j_i}=a)\no\\
&=\sum_{l=1}^{m-1}\sum_{\begin{subarray}{c}
                                    0\le i_1<...<i_l\le m-1\\
                                     1\le k_1,...,k_l\le n_0
                                  \end{subarray}}
   \sum_{(j_1,...,j_m)\in B(i_1,k_1,...,i_l,k_l)} P(Z_{j_1}=a)\no\\
   &\quad\quad\quad\quad\quad\quad\quad\quad\quad\quad\quad\quad\quad\quad\times
   \prod_{i=1}^{m-1}P(Z_{j_{i+1}}=a|Z_{j_i}=a).\label{gbd}
  \end{align}
 Fix $1\le l\le m,$ $0\le i_1<...<i_l\le m-1,$ and $1\le k_1,...,k_l\le n_0.$ If we can show
  \begin{align}
   \lim_{n\rto}\frac{1}{(\log n)^m} \sum_{(j_1,...,j_m)\in B(i_1,k_1,...,i_l,k_l)} P(Z_{j_1}=a)   \prod_{i=1}^{m-1}P(Z_{j_{i+1}}=a|Z_{j_i}=a)=0, \label{bik0}
  \end{align}
  then in view of \eqref{gbd} we can conclude that \eqref{gmn0} is true.

  In order to prove \eqref{bik0}, notice that
  \begin{align}
    &\sum_{(j_1,...,j_m)\in B(i_1,k_1,...,i_l,k_l)} P(Z_{j_1}=a)   \prod_{i=1}^{m-1}P(Z_{j_{i+1}}=a|Z_{j_i}=a)\no\\
    &\quad\quad=\sum_{(j_1,...,j_m)\in B(i_1,k_1,...,i_l,k_l)}\mathbf1_{i_1=0}P(Z_{j_1}=a) \prod_{s=2}^lP(Z_{j_{i_s+1}}=a|Z_{j_{i_s}}=a)\no\\
    &\quad\quad\quad\quad\quad\quad\quad\quad\quad\quad\quad\quad\times \prod_{i\in \{1,..,m-1\}/\{i_1,...,i_l\}}P(Z_{j_{i+1}}=a|Z_{j_i}=a)\no\\
    &\quad\quad +
    \sum_{(j_1,...,j_m)\in B(i_1,k_1,...,i_l,k_l)}\mathbf1_{i_1\ne 0}P(Z_{j_1}=a) \prod_{s=1}^lP(Z_{j_{i_s+1}}=a|Z_{j_{i_s}}=a)\no\\
    &\quad\quad\quad\quad\quad\quad\quad\quad\quad\quad\quad\quad\times \prod_{i\in \{1,..,m-1\}/\{i_1,...,i_l\}}P(Z_{j_{i+1}}=a|Z_{j_i}=a).\no
  \end{align}
From \eqref{fa}, we see that
\begin{align*}
  &\mathbf1_{i_1= 0}P(Z_{j_1}=a) \prod_{s=2}^lP(Z_{j_{i_s+1}}=a|Z_{j_{i_s}}=a)<c_2\mathbf1_{i_1= 0},\\
  &\mathbf1_{i_1\ne0} \prod_{s=1}^lP(Z_{j_{i_s+1}}=a|Z_{j_{i_s}}=a)<c_2\mathbf1_{i_1\ne 0},
\end{align*} for some constant $0<c_2<\infty$ independent of $n.$ Thus, on account of \eqref{pc} and \eqref{pj1} we get
\begin{align}
    &\sum_{(j_1,...,j_m)\in B(i_1,k_1,...,i_l,k_l)} P(Z_{j_1}=a)   \prod_{i=1}^{m-1}P(Z_{j_{i+1}}=a|Z_{j_i}=a)\no\\
    &\quad\quad\le\sum_{(j_1,...,j_m)\in B(i_1,k_1,...,i_l,k_l)}c_2\mathbf1_{i_1=0} \prod_{i\in \{1,..,m-1\}/\{i_1,...,i_l\}}\frac{1+\ve}{j_{i+1}-j_i}\no\\
    &\quad\quad +
    \sum_{(j_1,...,j_m)\in B(i_1,k_1,...,i_l,k_l)}c_2\mathbf1_{i_1\ne 0}\frac{1-\ve}{j_1}
    \prod_{i\in \{1,..,m-1\}/\{i_1,...,i_l\}}\frac{1+\ve}{j_{i+1}-j_i}\no\\
    &\quad\quad =c_2\sum_{1\le j_1'<...<j'_{m-l}\le n-\sum_{s=1}^lk_s}\prod_{i=0}^{m-l-1}\frac{(1-\ve)^{m-l}}{j'_{i+1}-j_i'}.\no
  \end{align}
Therefore, we conclude that
\begin{align}
   \varlimsup_{n\rto}&\frac{1}{(\log n)^m} \sum_{(j_1,...,j_m)\in B(i_1,k_1,...,i_l,k_l)} P(Z_{j_1}=a)   \prod_{i=1}^{m-1}P(Z_{j_{i+1}}=a|Z_{j_i}=a)\no\\
   &\le \varlimsup_{n\rto}\frac{ c_2(1-\ve)^{m-l}}{(\log n)^l}\frac{1 }{(\log n)^{m-l}}\sum_{1\le j_1'<...<j'_{m-l}\le n-\sum_{s=1}^lk_s}\frac{1}{j'_{i+1}-j_i'}=0\no
  \end{align}
where we use Lemma \ref{lemsa} to get the last step.
Thus, \eqref{bik0} is true and so is \eqref{gmn0}.

We turn next to prove \eqref{lgmn}. Dividing \eqref{gab} by $(\log n)^m$ and taking the lower and the upper limits, from \eqref{gmn1} and \eqref{gmn0}, we come to a conclusion that
\begin{align*}
  (1-\ve)^m<\varliminf \frac{G(m,n)}{(\log n)^m} \le \varlimsup\frac{G(m,n)}{(\log n)^m}<(1+\ve)^m.
\end{align*}
Since $\ve$ is arbitrary, letting $\ve\rightarrow 0,$ we get \eqref{lgmn}.

With the help of \eqref{lgmn}, from \eqref{eck} we deduce that
\begin{align}
  \lim_{n\rto}
  \frac{E|C\cap [1,n]|^k}{(\log n)^k}\lim_{n\rto}\sum_{m=1}^k\binom{k-1}{m-1}m!\frac{G(m,n)}{(\log n)^k}=k!\no
\end{align}
 which proves \eqref{cme}.

 We are now ready to prove \eqref{cd}. Using the fact $ ((2n)!)^{-\frac{1}{2n}}\sim \frac{e }{2n}$ as $n\rto,$ we get
\begin{align*}
  \sum_{k=1}^\infty ((2k)!)^{-\frac{1}{2k}}=\infty.
\end{align*}
Then, using Carleman's test for the uniqueness of the moment problem(see e.g. \cite[Chap.II, \S 12]{s}), we conclude that
\begin{align}
  \frac{|C\cap [1,n]|}{\log n}\rightarrow S\no
\end{align}
in distribution as $n\rto,$ where $S$ is an exponentially distributed random variable with $P(S>t)=e^{-t}, t>0.$ Theorem \ref{ce} is proved. \qed


\vspace{.5cm}

\noindent{{\bf \Large Acknowledgements:}}  I thank Professor Wenming Hong  for introducing to me the reference \cite{cfrb} from which I get inspiration to consider the present questions for BPVEs with immigration.

\end{document}